\documentclass{article}
\usepackage{amssymb}
\usepackage{amsmath}
\usepackage{amsfonts}
\usepackage{amscd}

\def\N{\mbox{I\hspace{-.15em}N}}

\def\R{\mbox{l\hspace{-.47em}R}}

\newcommand{\vas}{variables al\'{e}atoires }

\newtheorem{de}{\Large{D\'efinition}}[section]
\newtheorem{pp}{\Large{Proposition}}[section]

\newtheorem{lm}{\Large{Lemme}}[section]

\newtheorem{thm}{\Large{Th\'eor\`eme}}[section]
\newcommand{\dps}{\displaystyle }
\newcommand{\n}{\parallel }
\begin{document}
\begin{center}
{\bf\LARGE{Une g\'en\'eralisation de la notion d'int\'egrale it\'er\'ee relativement \`a un processus al\'eatoire.}}\\
Ludovic Valet\\
\small
{\textit{D\'epartement de Math\'ematiques, Universit\'e d'Angers, France.}}\\
\vspace{1,5cm}
{\bf Summary}\\
\end{center}
\footnotesize
\begin{itemize}
\item[]
\begin{itemize}
\item[]
In this paper we generalize notions of iterated integral with regard to an unpredictable process.
We establish a formula of integration by parts,  the existence of a continuous modification and give an expression of the increasing process.
\begin{center}
\hspace*{-2cm}{\bf Palabras llaves}

D\'ecomposition polyn\^omiale-Int\'egrale it\'er\'ee-Int\'egration par parties-Martingale-Modification continue-Produit sym´etrique-Produit tensoriel-Variation quadratique.
\end{center}
\begin{center}
\hspace*{-2cm}{\bf AMS}\\
60G12, 60H05, 60G50, 41A10, 33C45%
\end{center}
\end{itemize}
\end{itemize}
\normalsize
\section{Introduction}
Notre objectif est d'\'etendre le travail effectu\'e dans {{\bf [Schraf]}} et {{\bf [MFA\&ES]}}  en construisant puis \'etudiant une int\'egrale double.\\
Commençons par expliciter les notations de {{\bf [Schraf]}} pour ce qui nous concerne ici.
\begin{itemize}
\item
${\cal H}:=L^2(\R, B(\R), \mu)$ où $\mu$ est une mesure diffuse (pour le cas gaussien nous consid\`ererons $\mu=\lambda$ la mesure de Lebesgue). Nous noterons
$\langle,\rangle$ le produit scalaire et
$(e_n)_{n\in\N}$ une base hilbertienne.
\item
$(P_n)_{n\in\N}$ une famille de probabilit\'es sur $(\R, B(\R))$ telle que :
 $\int_{\R}xdP_n=0$,$\int_{\R}x^2dP_n=1$ et $\Bigl((\forall p>0)(\exists K_p)(\forall n\in\N)(\int_{\R}x^pdP_n\leq K_p)\Bigl)$.
Nous prendrons $\dps L^2(P):=L^2(\R^{\N},\dps\otimes_{{n\in\N}}B(\R),P)$ où $P=\otimes_{n\in\N}P_n$, noterons $X_n$ les
projecteurs canoniques de $\R^{\N}$ sur $\R$, $(n\in\N)$ et $H:=\overline{[X_n,\; n\in\N]}$.
Nous utiliserons la filtration $({\cal F}_n)_{n\in\N}:=(\sigma(X_k, k\leq n))_{n\in\N}$
\end{itemize}
Le mod\`ele suivant :
$$\begin{array}{ccccc}
\phi& : & {\cal H} & \to &H\\
&&h=\dps\sum_{k=0}^n\langle h, e_k\rangle e_k & \mapsto & \phi(h):=\dps\sum_{k=0}^n\langle h, e_k\rangle X_k\
\end {array}$$
va être \`a la base de notre \'etude. Pour des d\'etails de ses propri\'et\'es, on renvoie \`a {\bf [Schraf]}.\\
\\
Nous noterons $\Phi(h)_s:=\phi(h1_{]0,s]})$\\
Dans un premier temps, nous allons construire cette int\'egrale double. Pour ce faire il
nous faudra, au pr\'ealable, \'etudier dans $L^2(P)$, les produits $\phi(h)\phi(g)$. Nous terminerons
ce premier paragraphe par une g\'en\'eralisation du mod\`ele dans le but d'obtenir le d\'eveloppement
de ces produits dans une base de polyn\^omes orthogonaux et de faciliter le travail des
paragraphes suivants.\\
Les deuxi\`eme et troisi\`eme temps concernent l'\'etude du processus  $(Z_t)_{t\geq0}$ :
$$Z_t:=\int\Phi(h)_sd\Phi\left(g 1_{]0,t]}\right)_s=\int_0^t\Bigl(\int_0^sh(u)d\phi(u)\Bigl)g(v)d\phi(v)$$
que nous noterons plus simplement :
$Z_t:=\int_0^t\Phi(h)_sd\Phi\left(g \right)_s$.\\
(En cela nous ne suivrons pas exactement la chronologie de {\bf [Schraf]}, pr\'ef\'erant commencer par l'\'etablissement d'une formule d'int\'egration par partie et terminant par l'\'etablissement de l'existence d'une modification continue et de la variation quadratique de ce processus. Cette d\'emarche se  justifie techniquement.)
\section{Construction}
Nous voulons construire une int\'egrale double $\dps\int \Phi(h)_sd\Phi(g)_s$, et ainsi
g\'en\'eraliser la construction faite dans {\bf [Schraf]}. Cette construction s'appuie sur un
d\'evelopement adapt\'e dans une base de polyn\^omes des produits $\Phi(h)\Phi(g)$.\\
Cette int\'egrale sera prise dans le sens suivant :
$\dps\int \Phi(h)_sd\Phi(g)_s =\dps\lim_{L^2}\sum_{k=0}^{n}\Phi(h)_{t_{k}^n}\Bigl(\Phi(g)_{t_{k+1}^n}-\Phi(g)_{t_{k}^n}\Bigl)$
$=\dps\lim_{L^2}\sum_{k=0}^{n}\phi\Bigl(h1_{]0,t_{k}^n]}\Bigl)\phi\Bigl(g1_{]t_{k}^n,t_{k+1}^n]}\Bigl)$\\
Où $0=t_0^n<\cdots<t_n^n=t$ est une partition de $[0,t]$.\\
\\
Pour ce faire nous allons, dans la proposition 2.1., \'etudier les produits $\phi(h)\phi(g)$. Puis
montrer l'existence d'une limite, proposition 2.2., ce qui assurera l'existence de cette
int\'egrale comme un \'el\'ement de $L^2$.\\
Nous terminerons cette \'etude en remarquant qu'un changement de base de polyn\^omes
peut simplifier, dans certains cas, les calculs.
\subsection{Produits}
Pour un \'el\'ement f de ${\cal H \otimes  H}$, nous allons noter :
\begin{itemize}
\item[$\bullet$]
$\Biggl(
\dps\sum_{j=1}^{N}
\langle     f  ,    e_j   \otimes  e_j     \rangle   (X_j^2-1)
\Biggl)_{N\geq0}
=:\biggl(
\varphi_N^{(2)} (f)
\biggl)_{N\geq0}$
\item[$\bullet$]
$\Biggl(\dps\sum_{j=1}^{N}\Bigl(\sum_{k=0}^{j-1}\langle     f  ,    e_j   \otimes  e_k     \rangle  X_k   \Bigl)   X_j\Biggl)_{N\geq1}=:\biggl(\varphi_N^{(1,1)} (f)\biggl)_{N\geq1}$
\end{itemize}
Nous allons utiliser la base de polyn\^omes $\{(X_j^2-1)_{j\geq0};(X_jX_k)_{j\not=k}\}$.
D'une part $\biggl( \varphi_N^{(2)} (f) \biggl)_{N\geq0}$ et
$\biggl( \varphi_N^{(1,1)} (f)\biggl)_{N\geq1}$ sont des martingales relativement \`a la filtration
${\cal F}_j=\sigma(X_k; k\leq j)$.\\
On v\'erifie qu'il existe des limites ps et dans $L^2$ que nous noterons   $\varphi^{(2)} (f)$ et $\varphi^{(1,1)} (f)$.
D'autre part, des calculs simples d\'ecoulant des propri\'et\'es de martingales, montrent que  ces deux \'el\'ements sont orthogonaux dans $L^2$.
Le r\'esultat suivant, concernant les produits $\phi(h) \phi(g)$, sera \`a la base de la construction de l'int\'egrale double. Il d\'ecoule des propri\'et\'es mentionn\'ees pr\'ec\'edemment.
\begin{pp}
Pour des fonctions $f$ et $g$ de ${\cal H}$ on a :
\begin{enumerate}
\item
$\phi(h)\phi(g)\in L^2$
\item
$\phi(h)\phi(g) = \varphi^{(2)}(h\otimes g) + \varphi^{(1,1)}(h\otimes g+g\otimes h) +\langle h , g \rangle$ dans $L^2$
\end{enumerate} \end{pp}
\subsection*{Remarque}
On d\'etecte dans le point 2) le r\^ole particulier des fonctions sym\'etriques, qu'on retrouve
dans {\bf [PAM]} pour l'int\'egrale d'ordre 2.
\subsection{Limites}
Nous allons dans ce paragraphe, d\'eterminer la limite d'expressions du type :
$$\dps\sum_{k=0}^n\phi\Bigl(h1_{]0,t_{k}^n]}\Bigl)\phi\Bigl(g1_{]t_{k}^n,t_{k+1}^n]}\Bigl).$$
La proposition qui suit d\'ecoule des propri\'et\'es d'orthogonalit\'e, des propri\'et\'es de martingales et du th\'eor\`eme de convergence domin\'ee.

Pour expliciter la limite nous aurons besoin des deux indicatrices suivantes : 
$1_C =\{ (x,y)\in\R^2 \vert 0\leq x<y \leq t\}$\\
et\\
$1_{\widetilde{C}}=\{ (x,y)\in\R^2 \vert 0\leq y<x \leq t\}$.
\begin{pp}
Si nous prenons une suite de partitions de $[0,t]$ :
$$0=t^n_0<\cdots<t^n_{k_n}=t$$
dont le pas tend vers z\'ero, alors :
$$\dps \sum_{k=0}^{k_n}\phi \Bigl( h1_{]0,t^n_k]}\Bigl)\phi \Bigl( g1_{]t^n_k,t^n_{k+1}]}\Bigl)
\stackrel{L^2}{\to}
\varphi^{(1,1)}(h\otimes g 1_C +g\otimes h 1_{\widetilde{C}}) +\varphi^{(2)}
 (h\otimes g 1_C).$$
\end{pp}
\subsection*{\underline{\textbf{D\'emonstration}}}
\begin{enumerate}
\item[]
\small
$\dps\sum_{k=0}^{k_n}\phi\Bigl(h1_{  ]t^n_k,t^n_{k+1} ]} \Bigl) \phi \Bigl( g1_{]t^n_k,t^n_{k+1}]}\Bigl)=$\\
$\dps\sum_{k=0}^{k_n}
     \Bigl[
\varphi^{(1,1)}  \Bigl(    h\otimes g  1_{]0,t^n_k]}\otimes1_{  ]t^n_k,t^n_{k+1} ]} +
                                g\otimes h 1_{  ]t^n_k,t^n_{k+1} ]}\otimes1_{]0,t^n_k]}
                         \Bigl)$\\
 $\hspace*{1,5cm}+\varphi^{(2)}  \Bigl(    h\otimes g  1_{]0,t^n_k]}\otimes1_{  ]t^n_k,t^n_{k+1} ]}
                        \Bigl)
+\langle h1_{]0,t^n_k]}, g1_{  ]t^n_k,t^n_{k+1} ]}\rangle
\Bigl].$
\newline
On remarque que $\langle h1_{]0,t^n_k]}, g1_{  ]t^n_k,t^n_{k+1} ]}\rangle=0.$\newline
On va maintenant \'etudier les termes : \\
$\left\{
\begin{array}{lcl}
N^n_1   & := &
\left\Vert\dps\sum_{k=0}^{k_n}\varphi^{(1,1)}\biggl( h\otimes g  1_{]0,t^n_k]}\otimes1_{  ]t^n_k,t^n_{k+1} ]} \biggl)
-\varphi^{(1,1)}\biggl(h\otimes g1_C\biggl)\right\Vert_{L^2(\nu)}\\
N^n_2   & := &
\left\Vert
   \dps\sum_{k=0}^{k_n}
   \varphi^{(1,1)}\biggl(
                                g\otimes h  1_{  ]t^n_k,t^n_{k+1}]} \otimes1_{]0,t^n_k]}
                      \biggl)
 -  \varphi^{(1,1)}\biggl(
                                  g\otimes h1_{\widetilde{C}}
                       \biggl)
\right\Vert_{L^2(\nu)}\\
N^n_3   & := &
\left\Vert\dps\sum_{k=0}^{k_n}
\varphi^{(2)}\biggl(
                              h\otimes g  1_{]0,t^n_k]}\otimes1_{  ]t^n_k,t^n_{k+1} ]}
                    \biggl)
-\varphi^{(2)}\biggl(   h\otimes g1_C      \biggl)       \right\Vert_{L^2(\nu)}
\end{array}\right.$
\newline

$\begin{array}{lcl}
N^n_1   & = &
\left\Vert\dps\sum_{k=0}^{k_n}
\varphi^{(1,1)}\biggl(
                                 h\otimes g  1_{]0,t^n_k]}\otimes1_{  ]t^n_k,t^n_{k+1} ]} -   h\otimes g1_C
                      \biggl)
\right\Vert_{L^2(\nu)}\\
          & \leq &  \left\Vert h\otimes g
\biggl(
         \dps\sum_{k=0}^{k_n}
 1_{]0,t^n_k]}\otimes1_{  ]t^n_k,t^n_{k+1} ]} - 1_C
\biggl)
\right\Vert_{L^2(\mu^{\otimes^2})}\\
\end{array}$
\newline
Nous avons :
\begin{enumerate}
\item[$\bullet$]
$\biggl( h\otimes g\biggl)^2
\biggl(
          \dps\sum_{k=0}^{k_n} 1_{]0,t^n_k]}\otimes1_{  ]t^n_k,t^n_{k+1} ]} - 1_C
\biggl)^2 \leq
4\biggl( h\otimes g\biggl)^2\in L^2(\mu^{\otimes^2})$
\item[$\bullet$]
$\biggl(
          \dps\sum_{k=0}^{k_n} 1_{]0,t^n_k]}\otimes1_{  ]t^n_k,t^n_{k+1} ]} - 1_C
\biggl)^2\dps\mathop{\stackrel{\mu^{\otimes^2}}{\longrightarrow}}\limits_{\mbox{\ ps}}0$
\end{enumerate}
$\Longrightarrow
\biggl( h\otimes g\biggl)^2
\biggl(
          \dps\sum_{k=0}^{k_n} 1_{]0,t^n_k]}\otimes1_{  ]t^n_k,t^n_{k+1} ]} - 1_C
\biggl)^2\dps\mathop{\stackrel{\mu^{\otimes^2}}{\longrightarrow}}\limits_{\mbox{\ ps}}0$
\newline
En appliquant le th\'eor\`eme de convergence domin\'ee de Lebesgue on a : \\
$$
\left\Vert h\otimes g
\biggl(
         \dps\sum_{k=0}^{k_n}
 1_{]0,t^n_k]}\otimes1_{  ]t^n_k,t^n_{k+1} ]} - 1_C
\biggl)
\right\Vert_{L^2(\mu^{\otimes^2})}{\mathop{\longrightarrow}\limits_{n\to\infty}}0$$
d'où
$$\dps\sum_{k=0}^{k_n}\varphi^{(1,1)}\biggl( h\otimes g  1_{]0,t^n_k]}\otimes
1_{  ]t^n_k,t^n_{k+1} ]} \biggl)
\dps\mathop{  \stackrel{L^2(\nu)}{\longrightarrow} }\limits_{n\to\infty}
\varphi^{(1,1)}\biggl( h\otimes g1_C\biggl)
$$
Les termes $N_2^n$ et $N_3^n$ se traitent de la même façon.
\end{enumerate} \normalsize$\triangle$
Cette proposition nous permet de d\'efinir :%
\begin{de}
$$\int \Phi(h)_sd\Phi(g)_s:=\varphi^{(1,1)}\bigl(h\otimes g1_C+
g\otimes h1_{\widetilde{C}}\bigl)+\varphi^{(2)}\bigl(h\otimes g1_C\bigl)$$
\end{de}
Les sections qui suivent sont consacr\'ees \`a l'\'etude de cet objet. Avant de passer \`a cette
\'etude nous allons observer les cons\'equences d'un changement de base pour le d\'eveloppement des produits.
\subsection{Changement de base}
Nous allons d\'evelopper les s\'eries rencontr\'ees pr\'ec\'edemment dans la base des polyn\^omes
orthogonaux $(P_k)_k$ associ\'es aux \vas $(X_j)_j$. Cette \'ecriture nous permettra, en particulier,
 de mettre en \'evidence les propri\'et\'es concernant les fonctions sym\'etriques.\\
Nous allons, pour ce faire, g\'en\'eraliser le mod\`ele rappel\'e en introduction. Nous le ferons \`a un
 ordre quelconque tout d'abord, cela nous servira dans la suite, puis \`a l'ordre deux,
ce qui nous int\'eresse plus particuli\`erement dans ce paragraphe.\\
Dans l'\'etude faite dans {\bf [PAM]}, les polyn\^omes orthogonaux qui apparaissent naturellement
sont les polyn\^omes de Hermite (cas gaussien). Nous allons, pour retrouver les
propri\'et\'es du cas gaussien, prendre la diff\'erence des c\oe fficients des polyn\^omes issus
des $(X_j)_j$ (que nous noterons $\gamma_{..}$) et des polyn\^omes de Hermite
(que nous noterons $\Gamma_{..}$). Pour un compl\'ement sur les polyn\^omes orthogonaux on
 pourra consulter {\bf [Szego]}.\\
La construction de notre mod\`ele requiert l'introduction de deux op\'erateurs : \\
Fixons
$\left\{\begin{array}{lll}
\bullet & n\in\N^*\\
\bullet & (r,k)\in\{1,\ldots,n\}^2\\
\bullet & (\alpha_1,\ldots,\alpha_r)\in\N^r\mbox{\ \ tels que \ }\alpha_1+\cdots+\alpha_r=n\\
\bullet &(j_1,\ldots,j_r)\in\N^r\mbox{\ \ tels que \ }1\leq j_1<\cdots<j_r
\end{array}\right.$\\
On note ${\cal H}^{\circ^n}$ le produit tensoriel sym\'etrique et on d\'efinit les deux op\'erateurs
 qui suivent :
\begin{enumerate}
\item
$\Phi^{\circ^n} :\left\{ \begin{array}{ccl}
{\cal H}^{\circ^n} & \to          & L^2(P)\\
\mathop{\bigcirc}\limits_{i=1}^r e_{j_i}^{\circ^{\alpha_i}}
                 & \mapsto & P^{\alpha_1}(X_{j_1})\cdot\ldots\cdot P^{\alpha_r}(X_{j_r})
\end{array}\right.$\\
Où les $(P^{\alpha_i})_i$ sont les polyn\^omes associ\'es \`a la probabilit\'e $P_k$.
\item
$a_k^n :\left\{\begin{array}{lcl}
      {\cal H}^{\circ^n} & \to & H^{\circ^{n-k}}\\
\dps\mathop{\bigcirc}\limits_{i=1}^r e_{j_i}^{\circ^{\alpha_i}}  &\mapsto  &
\dps\sum_{\mathop{(k_1+\cdots+k_r=k)}\limits_{k_i\geq0}}
\prod_{(k_i\not=0)}\Bigl[\Bigl(
\gamma_{\alpha_i,\alpha_i-k_i} - \Gamma_{\alpha_i,\alpha_i-k_i}
\Bigl)1_{[\alpha_i\geq k_i]}\Bigl]
\dps\mathop{\bigcirc}\limits_{p=1}^re_{j_p}^{\alpha_p-k_p}
\end{array}\right.$
\end{enumerate}
On peut \'etendre ces deux op\'erateurs par lin\'earit\'e et continuit\'e \`a ${\cal H}^{\circ^n}$.\\
Ceux-ci nous permettent de donner l'\'ecriture suivante de notre int\'egrale double : \\
$$\dps\int \Phi(h)_sd\Phi(g)_s = \Bigl(
\Phi^{\circ^2}+\Phi\circ a_1^2
\Bigl)
\Bigl(\left[h\otimes g1_{C}\right]^{\circ}\Bigl)$$
Cette \'ecriture nous servira dans l'identification que l'on fera pour \'etudier l'existence d'une modification continue et le calcul de la variation quadratique.\\
\\
Dans le cas gaussien on retrouve la d\'ecomposition de {\bf [PAM]} :
$$\dps\int \Phi(h)_sd\Phi(g)_s =
\Phi^{\circ^2}
\Bigl(\left[h\otimes g1_{C}\right]^{\circ^2}\Bigl)$$
\subsection*{Remarques}
\begin{itemize}
\item
En passant sur l'espace $H^{\circ^n}$ nous n'avons plus une base normale. Cela va
nous obliger \`a transformer le produit scalaire pour conserver une isom\'etrie, ce que nous faisons ci-apr\`es.
\item
On peut \'egalement remarquer que pour n=1, $\Phi^{\circ^n}$ n'est
autre que $\Phi$.
\end{itemize}
\begin{itemize}
\item Nous prendrons le produit scalaire usuel sur $L^2(\mu)$ ainsi que
sur ${\cal H}^{\otimes^n}$.\vspace{0,5cm}
\item
On va consid\'erer le produit scalaire $\langle\, , \, \rangle_A^{{\cal H}^{\circ^n}}$
sur ${\cal H}^{\circ^n}$ d\'efini par :
$$\langle\, , \, \rangle_A^{{\cal H}^{\circ^n}}:=n!\langle A\cdot\, ,A\cdot\,
\rangle_{{\cal H}^{\otimes^n}}.$$
$\begin{array}{lccc}
\mbox{Où\, }A: &{\cal H}^{\circ^n} & \to      & {\cal H}^{\circ^n}\\
   & \mathop{\bigcirc}\limits_{i=1}^{r} e_{j_i}^{\alpha_i}
                & \mapsto  &
 \dps\sum_{\alpha} A_{ { j_1,\ldots,j_r}\atop{\alpha_1,\ldots,\alpha_r}}
 \mathop{\bigcirc}\limits_{i=1}^{r}e_{j_i}^{\alpha_i}
\end{array}$\\
avec
${A}_{ { j_1,\ldots,j_r}\atop{\alpha_1,\ldots,\alpha_r}}:=
\frac{\dps\prod_{i=1}^r E\left\vert P_{\alpha_i}(X_{j_i})\right\vert^2}
     {\dps\prod_{i=1}^r\alpha_i!}$.
\newline
Cet op\'erateur apparaît de façon naturelle dans la recherche d'une
isom\'etrie.
\end{itemize}
\subsection*{Notations}
\begin{itemize}
\item Posons
\begin{itemize}
\item
$C_{ {k_1,\ldots,k_r}\atop{ \alpha_1,\ldots,\alpha_r}}:=
\dps\prod_{k_i\not=0}\left[\left
(\gamma_{\alpha_i,\alpha_i-k_i}-\Gamma_{\alpha_i,\alpha_i-k_i}
\right)1_{[\alpha_i\geq k_i]}\right]$
\item
$C_{(k,n)}:=
\mathop{\sup}\limits_{k_1+\cdots+k_r=k\atop{\alpha_1+\cdots+\alpha_r=n}}
C_{ { \alpha_1,\ldots,\alpha_r}\atop{k_1,\ldots,k_r}}$
\item
$A_{(k,n)}:=
\mathop{\sup}\limits_{k_1+\cdots+k_r=k\atop{\alpha_1+\cdots+\alpha_r=n}}
A_{ {k_1,\ldots,k_r}\atop{ \alpha_1,\ldots,\alpha_r}}$
\end{itemize}
\item
On prend le produit scalaire usuel sur ${\cal H}^{\circ^n}$, \`a savoir :
$$<\;\cdot\; >_{{\cal H}^{\circ^n}}:=(n!)^2<\;\cdot\; >_{{\cal H}^{\otimes^n}}$$
\item
Les normes des op\'erateurs sont :
\begin{itemize}
\item
$\Vert T\Vert_{A}:= \mathop{\sup}\limits_{\Vert f\Vert\not=0}
\frac{\Vert T(f)\Vert_{L^2(\mu)}}{\Vert f \Vert_{A}}$
\item
$\Vert T_k\Vert_{{\cal H}^{\circ^n}}:=
\mathop{\sup}\limits_{\Vert f\Vert\not=0}
\frac{\Vert T_k(f)\Vert_{{\cal H}^{\circ{n-k}}}}{\Vert f\Vert_{{\cal H}^{\circ{n-k}}}}$
\end{itemize}
\end{itemize}
\section{Int\'egration par parties}
On cherche \`a \'etablir une formule du type It\^o pour les produits $\phi(h)\phi(g)$.\\
Celle-ci d\'ecoule naturellement de l'\'etude des produits faite en 2.1 et de la d\'efinition de
l'int\'egrale faite en 2.2.
\begin{thm}
Soient $h,g\in {\cal H}$, on a :
$$\Phi(h)\Phi(g)=\int\Phi(h)_sd\Phi(g)_s+\int\Phi(g)_sd\Phi(h)_s
+\langle h,g\rangle$$
\end{thm}
\subsection*{\underline{\textbf{D\'emonstration}}}
\begin{enumerate}
\item[]
\small
Posons $\Delta_{k,k+1}\Phi(h):=\Phi(h)_{t_{k+1}}-\Phi(h)_{t_{k}}$\\
$\dps\sum_{j=0}^{k_n} \Delta_{j,j+1}\Phi(h)
     \sum_{k=0}^{k_n}\Delta_{k,k+1}\Phi(g)
=   \dps\sum_{j=0}^{k_n}\Phi\left(h1_{]t_j^n,t_{j+1}^n]}\right)
     \sum_{k=0}^{k_n}\Phi\left(g1_{]t_k^n,t_{k+1}^n]}\right)\\
\begin{array}{cl}
= &\dps\sum_{j=0}^{k_n}\Phi\left(h1_{]t_j^n,t_{j+1}^n]}\right)
                   \Phi\left(g1_{]t_j^n,t_{j+1}^n]}\right)\\
&+  \dps\sum_{j=1}^{k_n}\sum_{k=0}^{j-1}\Phi\left(g1_{]t_k^n,t_{k+1}^n]}\right)
                   \Phi\left(h1_{]t_j^n,t_{j+1}^n]}\right)\\
&+  \dps\sum_{k=1}^{k_n}\sum_{j=0}^{k-1}\Phi\left(h1_{]t_j^n,t_{j+1}^n]}\right)
\Phi\left(g1_{]t_k^n,t_{k+1}^n]}\right).\\
\end{array}$
\newline
Appelons $S_i^n$ $(i\in\{1,2,3\})$ le $i^{\mbox{i\`eme}}$ terme du
membre de droite.
\begin{enumerate}
\item
On a :
$$S_1^n = \left(\Phi^{\circ^2}+\Phi\circ a_1^2\right)
\left(
h\otimes g \sum_{j=0}^{k_n} 1_{]t_j^n,t_{j+1}^n]}\otimes1_{]t_j^n,t_{j+1}^n]}
\right)$$
$$+\sum_{j=0}^{k_n} \langle
h1_{]t_j^n,t_{j+1}^n]},g 1_{]t_j^n,t_{j+1}^n]}
\rangle.$$
$\bullet$ D'une part :
$\dps\sum_{j=0}^{k_n}1_{]t_j^n,t_{j+1}^n]}\otimes1_{]t_j^n,t_{j+1}^n]}
\stackrel{L^2(\mu^{\otimes^2})}{\longrightarrow}1_{\Delta}$\\
\hphantom{0,5cm} et
$\langle
h\otimes g 1_{\Delta},e_{k_1}\otimes e_{k_2}
\rangle = 0.$\\
$\bullet$ D'autre part
$\dps\sum_{j=0}^{k_n}\langle
h1_{]t_j^n,t_{j+1}^n]},g1_{]t_j^n,t_{j+1}^n]}
\rangle = \langle h , g \rangle$.\\
Donc $S_1^n
\dps\mathop{\stackrel{L^2(\nu)}{\longrightarrow}}
\limits_{n\to\infty}
\langle h , g \rangle$
\item
$$S_2^n = \sum_{j=0}^{k_n}\phi\Bigl(g1_{  ]0,t_j]}\Bigl)\phi\Bigl(h1_{]t_j,t_{j+1}]}\Bigl)$$
En utilisant la propri\'et\'e 2.2. on a :
$$S_2^n
\dps\mathop{\stackrel{L^2(\nu)}{\longrightarrow}}
\limits_{n\to\infty}
\left(
\Phi^{\circ^2}+\Phi\circ a_1^2
\right)
\left(
g\otimes h1_C
\right)^{\circ^2}=:\int\Phi(g)_sd\Phi(h)_s.$$
\item
$$S_3^n =  \sum_{k=0}^{k_n}\phi\Bigl(h1_{  ]0,t^n_k]}\Bigl)\phi\Bigl(g1_{]t^n_k,t^n_{k+1}]}\Bigl)$$
de la même façon on a :
$$S_3^n
\dps\mathop{\stackrel{L^2(\nu)}{\longrightarrow}}
\limits_{n\to\infty}
\int\Phi(h)_sd\Phi(g)_s$$
\end{enumerate}
\end{enumerate} \normalsize$\triangle$\\
Toujours par les mêmes proc\'ed\'es et en \'ecrivant astucieusement le d\'eveloppement d'une
fonction dans ${\cal H}\otimes {\cal H}$, on peut calculer la norme de l'int\'egrale :
$$\left\Vert
\int\Phi(h)_sd\Phi(g)_s
\right\Vert_{L^2(\nu)}^2 =
\left\Vert
h\otimes g1_C
\right\Vert_{L^2(\nu)}^2
+\sum_{j\geq0}\langle h\otimes g,e_j\otimes e_j \rangle^2
\left(EX_j^4-3\right).$$
Ecriture qui fait imm\'ediatement apparaitre le r\'esultat pour l'int\'egrale de It\^o dans le cas
gaussien.
\section{R\'esultats pr\'eliminaires}
Pour \'etudier l'existence d'une modification continue et d\'eterminer la variation quadratique
nous aurons besoin d'un certain nombre d'identifications. C'est le propos de ce paragraphe.\\
\\
Plus pr\'ecis\'ement, nous voulons \'etudier le processus :
$$Z_t=\int_0^t\Phi(h)_sd\Phi\left(g \right)_s=(\Phi^{\circ^2}+
\Phi\circ a^2_1)(h_1\otimes(h_21_{]0,t]})1_C).$$
Pour ce faire nous aurons besoin de renseignements sur l'expression :
$$\mid Z_t-Z_s\mid^2.$$
Les r\'esultats qui suivent concernent l'\'etude de ce terme.
\subsection{Identification}
Le d\'eveloppement en s\'eries de $\mid Z_t-Z_s\mid^2$ est assez difficilement exploitable,
c'est pourquoi nous allons utiliser les op\'erateurs introduits
pr\'ec\'edemment et l'op\'erateur qui suit (cf {\bf [PAM]}) pour obtenir une \'ecriture de ce produit plus
exploitable :\\
$$\bigl(f\mathop{\sim}\limits_{1}f\bigl)(s_1,s_2):=\int Af(s_1,s_2)Af(s_3,s_2)d\mu(s_2)$$
Nous noterons $\pi_1$ la projection orthogonale sur $\{e_j\circ e_j \mid j\in \N\}$.\\
\\
Dans un premier temps on exprime $\mid Z_t-Z_s\mid^2$ dans la base des polyn\^omes
orthogonaux.\\
Puis on utilise les expressions de $f\circ f$, $a^n_k(f\circ f)$,
$f\mathop{\sim}\limits_{1}f$ et $a^n_k\bigl(f\mathop{\sim}\limits_{1}f\bigl)$
dans la base $(e_j\circ e_k)_{j,k}$ de ${\cal H}\circ {\cal H}$ pour identifier, ordre par ordre, chacun des
termes du d\'eveloppement pr\'ec\'edemment obtenu. On obtient la proposition suivante :
\begin{pp}
$$\left\vert\begin{array}{cll}\hline
&\bigl[(\Phi^{\circ^2}+\Phi\circ a_1^2)(f)\bigl]^2\\
 = & \Phi^{\circ^4}(f\circ f)&\mbox{(ordre 4)}\\
 +& \Phi^{\circ^3}\bigl(a_1^4(f\circ f)\Bigl) &\mbox{(ordre 3)}\\
 +&
\Phi^{\circ^2}\Bigl(4f\mathop{\sim}\limits_{1}f +a_2^4(f\circ f)
\Bigl) & \mbox{(ordre 2)}\\
 +& \Phi\bigl(
a_3^4(f\circ f) + 4 a_1^4(f\mathop{\sim}\limits_{1}f)-
6\bigl(a_1^4\circ \pi_1\bigl)(f\mathop{\sim}\limits_{1}f)
\Bigl) & \mbox{(ordre 1)}\\
+ & 2\n f\n^2+a_4^4(f\circ f) & \mbox{(ordre 0)}\\ \hline
\end{array}\right\vert$$
\end{pp}
\section{Modification continue}
Nous allons utiliser une g\'en\'eralisation du th\'eor\`eme de Kolmogorov-Centsov propos\'ee dans {\bf [Schraf]}, ainsi que l'identification \'etablie pr\'ec\'edemment.
\begin{lm}
$$E\left(\vert Z_t^n-Z_s^n\vert^2\right)^2
\leq(\frac{7}{2}C_{4,1}+C_{4,2}+C_{4,3}+C_{4,4}+2)
\n h_1\n_A^4\n h_21_{]s,t]}\n^4_A$$
\end{lm}
Pour la d\'emonstration, nous proc\'edons ordre par ordre. Les propri\'et\'es d'isom\'etrie
de $\Phi^{\circ^n}$ et de continuit\'e de $a_n^k$ et de $\pi_1$ nous fournissent,
en utilisant le d\'eveloppement de $(h_1\otimes h_21_{]s,t]}1_C)^{\circ^2}$  dans ${\cal H}\otimes {\cal H}$, des majorations rapides des termes :
$$K\mid\mid h_1\mid\mid^4\mid\mid h_21_{]s,t]}\mid\mid^4$$
où K est une constante.\\
\\
Ce lemme nous assure que nous sommes dans les conditions d'application du th\'eor\`eme de Kolmogorov-Centsov g\'en\'eralis\'e et l'on a le r\'esultat suivant :
\begin{thm}
Le processus $\left(\dps\int_0^t\Phi(h_1)d\Phi(h_2)\right)_{t}$
admet une modification continue.
\end{thm}
\section{Variation quadratique}
On veut montrer l'existence de la variation quadratique du processus :
$$\left(\dps\int_0^t\Phi(h_1)d\Phi(h_2)\right)_{t\geq0}.$$
Nous conservons la notation $Z_t:=\int_0^t\Phi(h_1)d\Phi(h_2)$.\\
Plus pr\'ecis\'ement, \'etant donn\'ee une suite de partitions
$0=t_0^n<t_1^n<\cdots<t_{k_n}^n=t$ de [0,t], nous cherchons la limite dans $L^2$
de la suite :
$$\left(\sum_{k=0}^{k_n}
\left\vert Z_{t^n_{k+1}}-Z_{t^n_k}
\right\vert^2
\right)_{n\in{\tiny\N}}$$
Nous avons d\'ej\`a \'etabli que (en notant $f_{k(n)}^{\circ}$ pour la sym\'etris\'ee de
$f_{t_k^nt_{k+1}^n}$)
$$\begin{array}{lcl}
\left\vert Z_{t^n_{k+1}}-Z_{t^n_k}
\right\vert^2 & = &
\Phi^{\circ^4}_n\bigl(f_{k(n)}^{\circ}\circ f_{k(n)}^{\circ}\bigl)\\
 & + &
\left(\Phi^{\circ^3}\circ a_1^2\right)_n
\bigl(f_{k(n)}^{\circ}\circ f_{k(n)}^{\circ}\bigl)\\
&+ &
\Phi^{\circ^2}_n\left(
4f_{k(n)}^{\circ}{\mathop{\sim}\limits_{1}}f_{k(n)}^{\circ}+
a^4_1\bigl(f_{k(n)}^{\circ}\circ f_{k(n)}^{\circ}\bigl)
\right)\\
&+&
\Phi_n\left(
a^4_3\bigl(f_{k(n)}^{\circ}\circ f_{k(n)}^{\circ}\bigl)+
4a^4_1\bigl(f_{k(n)}^{\circ}{\mathop{\sim}\limits_{1}}f_{k(n)}^{\circ}\bigl)
-6a^4_1\circ\pi_1\bigl(f_{k(n)}^{\circ}{\mathop{\sim}\limits_{1}}f_{k(n)}^{\circ}\bigl)
\right)\\
&+ &
2\n f_{k(n)}^{\circ}\n^2_A+
a^4_4\bigl(f_{k(n)}^{\circ}\circ f_{k(n)}^{\circ}\bigl)\\
\end{array}$$
La lin\'earit\'e et la continuit\'e des op\'erateurs  nous
incite, avant d'identifier la variation quadratique,
 \`a \'etudier, dans $L^2$,  les deux limites suivantes :
$$\left\{\begin{array}{l}
\dps\lim_{n\to\infty}\sum_{k=0}^{k_n}
f_{k(n)}^{\circ}\circ f_{k(n)}^{\circ}\\
\dps\lim_{n\to\infty}
4a^4_1\left(
\dps\sum_{k=0}^{k_n}\Bigl(f_{k(n)}^{\circ}{\mathop{\sim}\limits_{1}}f_{k(n)}^{\circ}
-6\pi_1\bigl(f_{k(n)}^{\circ}{\mathop{\sim}\limits_{1}}f_{k(n)}^{\circ}\bigl)
\Bigl)
\right)\\
\end{array}\right.$$
\begin{lm}
$$\lim_{n\to\infty}^{L^2}\dps\sum_{k=0}^{k_n}
f_{k(n)}^{\circ}\circ f_{k(n)}^{\circ}=0$$
\end{lm}
\subsection*{\underline{\textbf{D\'emonstration}}}
\begin{enumerate}
\item[]
\small
On a :\\
$f_{k(n)}^{\circ}\circ f_{k(n)}^{\circ}= \frac{1}{4}\Bigl(
h_1\otimes h^{k(n)}_2\otimes h_1\otimes h^{k(n)}_2
1_C\otimes1_C +
h_1\otimes h^{k(n)}_2\otimes h^{k(n)}_2\otimes h_1
1_C\otimes1_{\widetilde C}+$\\
$\hspace*{2,5cm}
h^{k(n)}_2\otimes h_1\otimes h_1\otimes h^{k(n)}_2
1_{\widetilde C}\otimes1_C +
h^{k(n)}_2\otimes h_1\otimes h^{k(n)}_2\otimes h_1
1_{\widetilde C}\otimes1_{\widetilde C}
\Bigl)$\\
$$\begin{array}{lcl}
\left\Vert
\dps\dps\sum_{k=0}^{k_n} f_{k(n)}^{\circ}\circ f_{k(n)}^{\circ}
\right\Vert^2  \leq
\dps\int
\left(
 \dps\sum_{k=0}^{k_n}\sum_{\sigma\in\Sigma_4}
 \Bigl(
  h_1\otimes h^{k(n)}_2\otimes h_1\otimes h^{k(n)}_2
  1_C\otimes1_C
 \Bigl)_{\sigma}
\right)^2d\mu^{\otimes^4}\\
\leq
(4!)^2\dps\sum_{\sigma\in\Sigma_4}
\Biggl[
 \dps\sum_{k=0}^{k_n}\int
 \Bigl(
  h_1\otimes h^{k(n)}_2\otimes h_1\otimes h^{k(n)}_2
  1_C\otimes1_C
 \Bigl)_{\sigma}^2d\mu^{\otimes^4}\\
+2\dps\sum_{k<j}\int
 \Bigl(
  h_1\otimes h^{k(n)}_2\otimes h_1\otimes h^{k(n)}_2
  1_C\otimes1_C
 \Bigl)_{\sigma}
 \Bigl(
  h_1\otimes h^{j(n)}_2\otimes h_1\otimes h^{j(n)}_2
  1_C\otimes1_C
 \Bigl)_{\sigma}d\mu^{\otimes^4}
\Biggl]\\
\hspace{5cm}\mbox{(in\'egalit\'e de convexit\'e)}\\
 \leq
(4!)^2\dps\sum_{\sigma\in\Sigma_4}
\Biggl[
 \dps\sum_{k=0}^{k_n}\int
 \Bigl[
  h_1^{\otimes^4}
 \Bigl]_{\sigma}d\mu^{\otimes^4}
 \int
 \Bigl[
  h_2^{\otimes^4}\Bigl(1_{]t_k^n,t_{k+1}^n]}\Bigl)^{\otimes^4}
 \Bigl]_{\sigma}d\mu^{\otimes^4}\\
+2\dps\sum_{k<j}\int
h_1^{\otimes^4}(x_{\sigma(1)},x_{\sigma(3)},y_{\sigma(1)},x_{\sigma(3)})
h_2^{\otimes^4}(x_{\sigma(2)},x_{\sigma(4)},y_{\sigma(2)},x_{\sigma(4)})\\
\hspace*{1cm}\cdot1_{C^{\times^4}}(x_{\sigma},y_{\sigma})
\Bigl(
 1_{]t_k^n,t_{k+1}^n]\times]t_j^n,t_{j+1}^n]}
 (x_{\sigma(2)},x_{\sigma(4)},y_{\sigma(2)},x_{\sigma(4)})
\Bigl)^{\otimes^2} \Biggl]d\mu^{\otimes^8}(x_{\sigma},y_{\sigma})
\end{array}$$
On sait que
$$\begin{array}{lcl}
\dps\dps\sum_{k=0}^{k_n}\Bigl(1_{]t_k^n,t_{k+1}^n]} \Bigl)^{\otimes^4}&
\mathop{\longrightarrow}\limits^{L^2(\mu^{\otimes^4})}&
1_{[x_{\sigma(2)}=x_{\sigma(4)}=y_{\sigma(2)}=y_{\sigma(4)}]}\\
\dps\sum_{k<j}
\Bigl(
 1_{]t_k^n,t_{k+1}^n]\times]t_j^n,t_{j+1}^n]}
\Bigl)^{\otimes^2}&
\mathop{\longrightarrow}\limits^{L^2(\mu^{\otimes^4})}&
1_{[x_{\sigma(2)}<x_{\sigma(4)},x_{\sigma(2)}=y_{\sigma(2)},
x_{\sigma(4)}=y_{\sigma(4)}]}\\
\end{array}$$
Pour \'etablir le dernier point il suffit de remarquer :\\
\begin{itemize}
\item d'une part que
$\left\{\begin{array}{lcl}
1_{[x=y]} & = &
\dps \lim_{n\to\infty}^{L^2}
\dps\sum_{k=0}^{k_n}1_{]t_k^n,t_{k+1}^n]\times]t_k^n,t_{k+1}^n]}(x,y)\\
1_{[x<y]} & = &
\dps \lim_{n\to\infty}^{L^2}
\sum_{k<j}1_{]t_k^n,t_{k+1}^n]\times]t_j^n,t_{j+1}^n]}(x,y)\\
\end{array}\right.$
\item et d'autre part que
$\left\{\begin{array}{l}
1_{[x=y]} 1_{[z=t]}1_{[x<z]}1_{[y<t]} \\
= \dps \lim_{n\to\infty}^{L^2}\Biggl[
\dps\sum_{k=0}^{k_n}1_{]t_k^n,t_{k+1}^n]^2}(x,y)\cdot
\dps\sum_{k=0}^{k_n}1_{]t_k^n,t_{k+1}^n]^2}(z,t)\cdot\\
\dps\sum_{k<j}1_{]t_k^n,t_{k+1}^n]\times]t_j^n,t_{j+1}^n]}(x,z)\cdot
\sum_{k<j}1_{]t_k^n,t_{k+1}^n]\times]t_j^n,t_{j+1}^n]}(y,t)\Biggl]
\\
\end{array}\right.$
\item car
$\left\{\begin{array}{l}
\Biggl[\dps\dps\sum_{k=0}^{k_n}1_{]t_k^n,t_{k+1}^n]^2}(x,y)\cdot
\dps\sum_{k=0}^{k_n}1_{]t_k^n,t_{k+1}^n]^2}(z,t)\cdot\\
\dps\sum_{k<j}1_{]t_k^n,t_{k+1}^n]\times]t_j^n,t_{j+1}^n]}(x,z)\cdot
\sum_{k<j}1_{]t_k^n,t_{k+1}^n]\times]t_j^n,t_{j+1}^n]}(y,t)\Biggl]\\
=\dps\sum_{k<j}1_{]t_k^n,t_{k+1}^n]^2}\otimes1_{]t_j^n,t_{j+1}^n]^2}(x,y,z,t)\\
\end{array}\right.$
\end{itemize}
Finalement, en utilisant le th\'eor\`eme de convergence domin\'ee de
Lebesgue et la mesure nulle de ces ensembles on a le r\'esultat
annonc\'e.
\end{enumerate} \normalsize$\triangle$\\
Soient $h_1,\;h_2\in{\cal H}$ et $0=t_0^n<t_1^n<\cdots<t_{k_n}^n=t$ une suite de
partitions de $[0,t]$ dont le pas tend vers 0.\\
Nous noterons $f_{k(n)}:=h_1\otimes \left(h_21_{]0,t_k^n]}\right)\cdot1_C$
\begin{lm}
Pour $\Psi=\Phi^{\circ^2}$ ou $\Psi=\Phi\circ a_1^4$ on a :
$$\mathop{\lim}\limits_{n\to\infty}E\left[
4\Psi
\bigl(\dps\sum_{k=0}^{k_n}f_{k(n)}\mathop{\sim}\limits_{1}f_{k(n)}\bigl)-
\int_0^th_2^2\Psi
\Bigl(\bigl[h_11_{]0,.]}
\bigl]^{\circ^2}\Bigl)d\mu\right]^2=0$$
\end{lm}
\subsection*{\underline{\textbf{D\'emonstration}}}
\begin{enumerate}
\item[]
\small
Nous regarderons le cas $\Psi=\Phi^{\circ^2}$. Nous
Introduisons quelques notations :\\
\begin{itemize}
\item
$S_n:=\left[
4\Phi^{\circ^2}
\bigl(\dps\sum_{k=0}^{k_n}f_{k(n)}\mathop{\sim}\limits_{1}f_{k(n)}\bigl)-
\int_0^th_2^2\Phi^{\circ^2}
\Bigl(\bigl[h_11_{]0,.]}
\bigl]^{\circ^2}\Bigl)d\mu\right]^2$
\item
$h^{k(n)}:=h1_{]0,t_k^n]}$
\end{itemize}
$$\begin{array}{clc}
&4f_{k(n)}\mathop{\sim}\limits_{1}f_{k(n)}(s_1,s_3)=\\
&\int h_1^{\otimes^2}(s_1,s_2)h_2^{{k(n)}\otimes^2}(s_2,s_3)
1_C(s_1,s_2)1_C(s_2,s_3)d\mu(s_2)  &
\bigl(:=I_1^{k(n)}(s_1,s_2)\bigl)\\
+&\int h_1^{\otimes^2}(s_1,s_3)h_2^{{k(n)}\otimes^2}(s_2,s_2)
1_C(s_1,s_2)1_C(s_3,s_2)d\mu(s_2)  &
\bigl(:=I_2^{k(n)}(s_1,s_2)\bigl)\\
+&\int h_1^{\otimes^2}(s_2,s_2)h_2^{{k(n)}\otimes^2}(s_1,s_3)
1_C(s_2,s_1)1_C(s_2,s_3)d\mu(s_2)  &
\bigl(:=I_3^{k(n)}(s_1,s_2)\bigl)\\
+&\int h_1^{\otimes^2}(s_2,s_3)h_2^{{k(n)}\otimes^2}(s_1,s_2)
1_C(s_2,s_1)1_C(s_3,s_2)d\mu(s_2)  &
\bigl(:=I_4^{k(n)}(s_1,s_2)\bigl)\\
\end{array}$$
L'in\'egalit\'e de convexit\'e nous permet d'\'ecrire :
$$S_n\leq
4E\Bigl[\Phi^{\circ^2}\dps\sum_{k=0}^{k_n}I^{k(n)}_1\Bigl]^2
+4E\Bigl[\Phi^{\circ^2}\dps\sum_{k=0}^{k_n}I^{k(n)}_2-
\int_0^th^2_2\Phi^{\circ^2}\bigl[(h_11_{]0,.]})^{\circ^2}\bigl]d\mu
\Bigl]^2$$
$$+4E\Bigl[\Phi^{\circ^2}\dps\sum_{k=0}^{k_n}I^{k(n)}_3\Bigl]^2
+4E\Bigl[\Phi^{\circ^2}\dps\sum_{k=0}^{k_n}I^{k(n)}_4\Bigl]^2
$$
\fbox{Etape 1 : } Calcul de
$\mathop{\lim}\limits_{n\to\infty}E\Bigl[\Phi^{\circ^2}
\dps\sum_{k=0}^{k_n}I^{k(n)}_1\Bigl]^2$\\
$$\begin{array}{cl}
&E\Phi^{\circ^2}(I_1^{k(n)})^2 =\n
I_1^{k(n)}\n^2_{L^2(\mu^{\otimes^2})}\hspace{2cm} \mbox{(isom\'etrie)}\\
= & \int\Bigl(\dps\dps\sum_{k=0}^{k_n}
\int h_1^{\otimes^2}(s_1,s_2)h_2^{{k(n)}\otimes^2}(s_2,s_3)
1_C(s_1,s_2)1_C(s_2,s_3)d\mu(s_2)
\Bigl)^2d\mu^{\otimes^2}(s_1,s_3)\\
= &\int\Bigl[\dps\dps\sum_{k=0}^{k_n}
\int h_1^{\otimes^4}(s_1,s_2,s_1,t_2)h_2^{{k(n)}\otimes^4}(s_2,s_3,t_2,s_3)
1_{C^4}(s_1,s_2,s_2,s_3,s_1,t_2,t_2,s_3)d\mu^{\otimes^2}(s_2,t_2)
\\
& +2\dps\sum_{k(n)<j(n)}
\int h_1^{\otimes^4}(s_1,s_2,s_1,t_2)
h_2^{k(n)\otimes^2}(s_2,s_3)h_2^{j(n)\otimes^2}(t_2,s_3)\\
&\hspace{1,5cm}1_{C^4}(s_1,s_2,s_2,s_3,s_1,t_2,t_2,s_3)d\mu^{\otimes^2}(s_2,t_2)
\Bigl]d\mu^{\otimes^2}(s_1,s_3)\\
\end{array}$$
\begin{itemize}
\item
D'une part $\dps\dps\sum_{k=0}^{k_n}1_{]t_k^n,t_{k+1}^n]^4}(s_2,s_3,t_2,s_3)
\mathop{\longrightarrow}\limits_{n\to\infty}^{L^2(\mu^{\otimes^4})}
1_{[s_2=s_3=t_2]}$
\item
D'autre part :
$\dps\sum_{k(n)<j(n)}1_{]t_k^n,t_{k+1}^n]^2\times[t_j^n,t_{j+1}^n]^2}(s_2,s_3,t_2,s_3)=0$\\
car $(k(n)<j(n))\Rightarrow(]t_k^n,t_{k+1}^n]\bigcap[t_j^n,t_{j+1}^n]=\emptyset)$
\end{itemize}
Le th\'eor\`eme de convergence domin\'ee de Lebesgue nous donne :
$$\mathop{\lim}\limits_{n\to\infty}E\Bigl[\Phi^{\circ^2}
\dps\sum_{k=0}^{k_n}I^{k(n)}_1\Bigl]^2=\int h_1^{\otimes^2}(s_1,s_2)h_2^{{k(n)}\otimes^2}(s_2,s_3)
1_{[s_2=s_3=t_2]}
d\mu^{\otimes^4}(s_1,s_2,t_2,s_3)=0$$
\fbox{Etape 2} Calcul de
$\mathop{\lim}\limits_{n\to\infty}E\Bigl[\Phi^{\circ^2}
\dps\sum_{k=0}^{k_n}I^{k(n)}_i\Bigl]^2$ pour $i\in\{3,4\}$\\
On proc\`ede de la même façon et on a \'egalement :
$$\mathop{\lim}\limits_{n\to\infty}E\Bigl[\Phi^{\circ^2}
\dps\sum_{k=0}^{k_n}I^{k(n)}_3\Bigl]^2=\mathop{\lim}\limits_{n\to\infty}E\Bigl[\Phi^{\circ^2}
\dps\sum_{k=0}^{k_n}I^{k(n)}_4\Bigl]^2=0$$
\fbox{Etape 3 : } calcul de $\mathop{\lim}\limits_{n\to\infty}
E\Bigl[\Phi^{\circ^2}\dps\sum_{k=0}^{k_n}I^{k(n)}_2-
\int_0^th^2_2\Phi^{\circ^2}\bigl[(h_11_{]0,.]})^{\circ^2}\bigl]d\mu
\Bigl]^2$\\
$$\begin{array}{lcl}
\dps\dps\sum_{k=0}^{k_n}I_2^{k(n)}(s_2,s_3) & = &
\dps\dps\sum_{k=0}^{k_n}\int_{t_k^n}^{t_{k+1}^n}
h_2^2(s_2)h_1^{\otimes^2}(s_1,s_3)1_C(s_2,s_3)1_{\widetilde C}(s_2,s_3)
d\mu(s_2)\\
& = &\dps\int_0^th_2^2(s_2)h_1^{\otimes^2}(s_1,s_3)
1_C(s_2,s_3)1_{\widetilde C}(s_2,s_3)
d\mu(s_2)\\
& = &\dps\int_0^th_2^2(s_2)
\Bigl(h_11_{]0,s_2]}\otimes h_11_{]0,s_2]}
\Bigl)(s_1,s_3)d\mu(s_2)\\
\end{array}$$
Il s'agit donc de montrer que :
$$\Phi^{\circ^2}\Bigl(
\int_0^th_2^2(s_2)
h_11_{]0,s_2]}\otimes h_11_{]0,s_2]}d\mu(s_2)\Bigl)
\stackrel{L^2}{=}
\int_0^th_2^2(s_2)
\Phi^{\circ^2}\Bigl(\bigl(
h_11_{]0,s_2]}\bigl)^{\circ^2}\Bigl)d\mu(s_2)
$$
Pour cela on va d\'ecomposer $1_C$ en produit tensoriel :
$$1_C (s_1,s_2) \stackrel{L^2}{=}
\mathop{\lim}\limits_{n\to\infty}\sum_{k=0}^{k_n}
1_{]0,t_{k-1}^n]\times]t_k^n,t_{k+1}^n]}(s_1,s_2) $$
$$1_{\widetilde C} (s_3,s_2) \stackrel{L^2}{=}
\mathop{\lim}\limits_{n\to\infty}\sum_{k=0}^{k_n}1_{]t_k^n,t_{k+1}^n]\times]0,t_{k-1}^n]}(s_3,s_2) $$
Nous allons noter $C_1^{k(n)}:=1_{]0,t_{k-1}^n]}$ et $C_2^{k(n)}:=1_{]t_k^n,t_{k+1}^n]}$
$$\begin{array}{cl}
&\int_0^th_2^2(s_2)h_1^{\otimes^2}(s_1,s_3)
1_C(s_2,s_3)1_{\widetilde C}(s_2,s_3)
d\mu(s_2)\\
= &
h_1^{\otimes^2}\int_0^th_2^2(s_2)
\mathop{\lim}\limits_{{n\to\infty}}\dps\sum_{k=0}^{k_n}\sum_{l=0}^{k_n}
1_{C_1^{k(n)}\times C_2^{k(n)}}(.,s_2)
1_{C_2^{l(n)}\times C_1^{l(n)}}(s_2,.)d\mu(s_2)\\
& \mathop{=}\limits_{{\mbox{\scriptsize cvg}}\atop{\mbox{\tiny domin\'ee}}}
\mathop{\lim}\limits_{n\to\infty}\dps\sum_{{k(n)},{l(n)}}
h_1^{\otimes^2}1_{C_1^{k(n)}}1_{C_1^{l(n)}}
\int_0^th_2^2(s_2)1_{C_2^{k(n)}}(s_2)1_{C_2^{l(n)}}(s_2)d\mu(s_2)
\end{array}$$
$$\begin{array}{clc}
&\Phi^{\circ^2}\Bigl(\dps
\int_0^th_2^2(s_2)h_1^{\otimes^2}1_C(.,s_2)1_{\widetilde
C}(s_2,.)d\mu(s_2)\Bigl) \\
=&\Phi^{\circ^2}\left(
\mathop{\lim}\limits_{n\to\infty}\dps\sum_{k(n),l(n)}
h_1^{\otimes^2}1_{C_1^{k(n)}}1_{C_1^{l(n)}}
\int_0^th_2^2(s_2)1_{C_2^{k(n)}}(s_2)1_{C_2^{l(n)}}(s_2)d\mu(s_2)\right) \\
= &\mathop{\lim}\limits_{n\to\infty}
\Phi^{\circ^2}
\left(
\dps\sum_{k(n),l(n)}h_1^{\otimes^2}1_{C_1^{k(n)}}1_{C_1^{l(n)}}
\int_0^th_2^2(s_2)1_{C_2^{k(n)}}(s_2)1_{C_2^{l(n)}}(s_2)d\mu(s_2)
\right) &\\
       &\hspace*{3cm} \left\{\begin{array}{l}
        t\mbox{ fix\'e}\\
        \mbox{cont. de } \Phi^{\circ^2}\\ \mbox{(Prop. 2.1.4.)}
        \end{array}\right.\\
= &\mathop{\lim}\limits_{ n\to\infty}
\Phi^{\circ^2}
\left(
\dps\sum_{k(n),l(n)}h_1^{\otimes^2}1_{C_1^{k(n)}}1_{C_1^{l(n)}}
\right)
\int_0^th_2^2(s_2)1_{C_2^{k(n)}}(s_2)1_{C_2^{l(n)}}(s_2)d\mu(s_2)\\
=&\mathop{\lim}\limits_{n\to\infty}
\dps\int_0^t
\Phi^{\circ^2}
\left(
\dps\sum_{k(n),l(n)}h_1^{\otimes^2}1_{C_1^{k(n)}}1_{C_1^{l(n)}}
\right)
h_2^2(s_2)1_{C_2^{l(n)}}(s_2)1_{C_2^{l(n)}}(s_2)d\mu(s_2)\\
\end{array}$$
On peut donc \'ecrire :
$$E\Bigl[\dps\sum_{k=0}^{k_n}\Phi^{\circ^2}(I_2^{k(n)})-
\int_0^th_2^2(s_2)\Phi^{\circ^2}
\bigl[\bigl(h_11_{]0,s_2]}\bigl)^{\circ^2}\bigl]d\mu(s_2)
\Bigl]^2$$
$$=E\Biggl[\mathop{\lim}\limits_{n\to\infty}
  \Bigl[\int_0^th_2^2(s_2)
   \Bigl(\Phi^{\circ^2}
   \bigl(
   h_1\otimes h_1 \sum_{(k(n),l(n))}
   1_{C_1^{k(n)}\times C_2^{k(n)}}(.,s_2)1_{C_2^{l(n)}\times C_1^{l(n)}}(s_2,.)
   \bigl)-$$
   $$\Phi^{\circ^2}
   \bigl(
   h_1\otimes h_1
   1_{C}(.,s_2)1_{\widetilde C}(s_2,.)
   \bigl)
   \Bigl)d\mu(s_2)
  \Bigl]
\Biggl]^2$$
Dans la suite on va noter :
\begin{itemize}
\item
$Z_{k(n)}(.,.,s_2) :=\Phi^{\circ^2}
   \bigl(
   h_1\otimes h_1\dps \sum_{(k(n),l(n))}
   1_{C_1^{k(n)}\times C_2^{k(n)}}(.,s_2)1_{C_2^{l(n)}\times C_1^{l(n)}}(s_2,.)
   \bigl)$
\item
$Z(.,.,s_2):=\Phi^{\circ^2}
   \bigl(
   h_1\otimes h_1
   1_{C}(.,s_2)1_{\widetilde C}(s_2,.)
   \bigl)$
\item
$Y_{k(n)} := \vert Z_{k(n)}-Z\vert$
\end{itemize}
$$\begin{array}{cl}
& E\Biggl[
  \Bigl[\dps\int_0^th_2^2(s_2)
   \Bigl(\Phi^{\circ^2}
   \bigl(
   h_1\otimes h_1 \sum_{(k(n),l(n))}
   1_{C_1^{k(n)}\times C_2^{k(n)}}(.,s_2)1_{C_2^{l(n)}\times C_1^{l(n)}}(s_2,.)
   \bigl)-\\
   &\hspace{1cm}\Phi^{\circ^2}
   \bigl(
   h_1\otimes h_1
   1_{C}(.,s_2)1_{\widetilde C}(s_2,.)
   \bigl)
   \Bigl)d\mu(s_2)
  \Bigl]
\Biggl]^2\\
\leq & E\left( \dps\int_0^t h_2^2(s_2)
\left\vert
Z_{k(n)}(.,.,s_2)-Z(.,.,s_2)
\right\vert d\mu(s_2)\right)^2\\
=&E\left( \dps\int_0^t\int_0^t \bigl(h_2^2\bigl)^{\otimes^2}(s_2,t_2)
Y_{k(n)}(s_2)Y_{k(n)}(t_2)d\mu(s_2,t_2)\right)\\
=&\dps\int_0^t\int_0^t \bigl(h_2^2\bigl)^{\otimes^2}(s_2,t_2)
E\left( Y_{k(n)}(s_2)Y_{k(n)}(t_2)\right)d\mu(s_2,t_2)\\
\leq &\dps\int_0^t\int_0^t \bigl(h_2^2\bigl)^{\otimes^2}(s_2,t_2)
\left(E \left\vert Y_{k(n)}(s_2)\right\vert^2\right)^{\frac{1}{2}}
\left(E \left\vert Y_{k(n)}(t_2)\right\vert^2\right)^{\frac{1}{2}}
d\mu(s_2,t_2)\\
\end{array}$$
Pour finir il suffit de voir que :\\
$
E \left\vert Y_{k(n)}(s_2)\right\vert^2  =
E\Bigl[
\Phi^{\circ^2}\bigl(h_1\otimes h_1\dps\sum_{(k(n),l(n))}
1_{C_1^{k(n)}\times C_2^{k(n)}}(.,s_2)1_{C_2^{l(n)}\times C_1^{l(n)}}(s_2,.)-$\\
$\hspace*{4cm}h_1\otimes h_1\dps\sum_{(k(n),l(n))}
1_{C}(.,s_2)1_{\widetilde C}(s_2,.)
\bigl)\Bigl]^2$\\
 $= \left\Vert h_1\otimes h_1
 \Bigl[
 1_{C}(.,s_2)1_{\widetilde C}(s_2,.)-
\dps\sum_{(k(n),l(n))}
1_{C_1^{k(n)}\times C_2^{k(n)}}(.,s_2)1_{C_2^{l(n)}\times C_1^{l(n)}}(s_2,.)
\Bigl]
 \right\Vert^2
$
\end{enumerate} \normalsize$\triangle$\\
Par ailleurs de simples calculs donnent : \\
$$\Biggl(\Phi\Bigl(h_11_{]0,s]}\Bigl)\Biggl)^2=
\Phi^{\circ^2}\Bigl(\Bigl[h_11_{]0,s]}\Bigl]^{\circ^2}\Bigl)+
\Phi\circ a_1^2\Bigl(\Bigl[h_11_{]0,s]}\Bigl]^{\circ^2}\Bigl)+
\left\Vert
h_11_{]0,s]}
\right\Vert^2$$
Ce qui a l'aide du paragraphe pr\'eliminaire nous permet d'\'enoncer :
\begin{thm}
$$\lim_{n\to\infty}^{L^2}\dps\sum_{k=0}^n\left\vert Z_{t_{k+1}^n}-Z_{t_{k}^n}\right\vert^2
= \int_0^th_2^2(s)\Bigl(\Phi\bigl(h_11_{[0,s]}\bigl)\Bigl)^2d\mu(s)+
\int_0^th_2^2(s)
\Phi\circ a_1^2\circ \pi_1
\bigl(  h_11_{]0,s]}  \bigl)^{\circ^2}d\mu(s).$$
\end{thm}
\subsection*{Remarque}
\begin{itemize}
\item
Dans le cas où les variables $\bigl(X_k\bigl)_k$ suivent une loi normale
on retrouve la variation quadratique gaussienne puisque dans ce
cas $a_1^2\circ\pi_1\bigl(h_11_{]0,s]}\bigl)^{\circ^2}=0.$
\end{itemize}
\newpage
\begin{center}
{\LARGE BIBLIOGRAPHIE}
\end{center}
\vspace{3cm}
{\bf {\bf [Szego]}}   SZEGO :   \textit{Orthogonal polynomial}, 1939
\newline
{\bf {\bf [Schraf]}}    SCHRAFSTETTER, Eric :   \textit{Quelques aspects des s\'eries al\'eatoires
d\'efinies par une mesure vectorielle, application aux EDPS}, Th\`ese
de doctorat, Angers, 1998.
\newline
{{\bf [MFA\&ES]}} ALLAIN, Marie-France, SCHRAFSTETTER, Eric : \textit{Quelques aspects des s\'eries al\'eatoires d\'efinies par une mesure
 vectorielle}, Pr\'epublication du D\'epartement de Math\'ematiques, Universit\'e d'Angers, $\mbox{n}^{\circ}41$ (1997).
\newline
{\bf {\bf [PAM]}} MEYER, Paul Andr\'e    :   \textit{Quantum Probability for
Probabilists}, Springer, LNM 1533, 1995.
\end{document}